\documentclass[12pt]{article}
\usepackage[left=0.65in, top = 0.8in, bottom = 0.8in, right = 0.65in]{geometry}

\usepackage{adjustbox}
\usepackage{enumitem}
\usepackage{multirow}
\usepackage{pdflscape}
\usepackage{booktabs}
\usepackage{floatrow}
\usepackage[counterclockwise, figuresright]{rotating}
\usepackage{enumitem}
\usepackage{hyperref}
\usepackage{subcaption}
\usepackage{eurosym}
\usepackage{natbib}
\setcitestyle{numbers,square}
\usepackage{amsmath,amsfonts,amssymb}
\DeclareUnicodeCharacter{20AC}{\EUR}
\usepackage[margin=1in]{caption}
\usepackage{float}
\usepackage{rotating}
\usepackage{arydshln}
\usepackage{longtable}
\usepackage{tabulary}
\usepackage{tabularx}
\usepackage{chngcntr}

\usepackage{setspace}
\usepackage{graphicx}
\usepackage{bm}
\usepackage{amsmath}
\usepackage{adjustbox}

\usepackage{amsmath}
\usepackage{amssymb}
\usepackage{amsthm}
\usepackage{amscd}
\usepackage{amsfonts}
\usepackage{eurosym}
\usepackage{bezier}

\usepackage{enumerate}
\usepackage{euro}
\usepackage{epsfig}
\usepackage{layout}

\usepackage{natbib}

\usepackage{float}
\usepackage{multirow}
\usepackage[flushleft]{threeparttable}
\usepackage{booktabs}
\usepackage[marginal]{footmisc}
\usepackage{lscape}
\usepackage{makecell}
\usepackage{longtable}
\usepackage{tabulary}
\usepackage{tabularx}
\newcolumntype{Y}{>{\raggedleft\arraybackslash}X}

\makeatletter
\renewcommand{\fnum@figure}{Fig. \thefigure}
\makeatother

\usepackage{xcolor}
\hypersetup{
	colorlinks,
	linkcolor={blue!50!black},
	citecolor={blue!50!black},
	urlcolor={blue!80!black},
	breaklinks=true
} % define hyperlink colors

\usepackage[nameinlink,capitalise]{cleveref}
\Crefname{figure}{Fig.}{Figs.}
\Crefname{equation}{Eq.}{Eqs.}

\usepackage{lipsum}
\usepackage{titlesec}
\usepackage{tcolorbox}
\tcbuselibrary{skins}

\definecolor{almond}{rgb}{0.94, 0.87, 0.8}
\definecolor{brass}{rgb}{0.71, 0.65, 0.26}
\definecolor{gold(metallic)}{rgb}{0.83, 0.69, 0.22}
\usepackage{authblk}

\begin{document}

\title{On Finite Time Span Estimators of Parameters for Ornstein-Uhlenbeck Processes}
\author{Jun S. Han, Nino Kordzakhia \\ \small School of Mathematical and Physical Sciences, Macquarie University \\
\small \textit{(email: jun.han1@hdr.mq.edu.au, nino.kordzakhia@mq.edu.au)}}

\date{}
\maketitle

\begin{abstract}
    We study the bias and the mean-squared error of the maximum likelihood estimators (MLE) of parameters associated with a two-parameter mean-reverting process for a finite time $T$. Using the likelihood ratio process, we derive the expressions for MLEs, then compute the bias and the MSE via the change of measure and Ito's formula. We apply the derived expressions to the general Ornstein-Uhlenbeck process, where the bias and the MSE are numerically computed through a joint moment-generating function of key functionals of the O–U process. A numerical study is provided to illustrate the behaviour of bias and the MSE for the MLE of the mean-reverting speed parameter. 
\end{abstract}

\section{Introduction}

We focus on the univariate two-parameter mean-reverting process $\mathbf{Y}=\{Y_{t},0\leq t\leq T\}$, which is defined through the following stochastic differential equation (SDE),
\begin{equation}
dY_{t}=\lambda (\alpha -b_{t})dt+\sigma_{t}dW_{t},\quad Y_{0}=y,\quad\lambda \in {\Lambda},\quad \alpha \in \mathcal{A}.  \label{Mainmodel}
\end{equation}
The two parameters $\lambda$ and $\alpha$ are to be estimated from the following sets of observations $\{Y_{t},b_{t},\sigma _{t},\; 0\leq t\leq T\} $ in open intervals $\Lambda\in\mathbb{R}$ and $\mathcal{A}\in\mathbb{R}$, respectively. The standard Brownian motion $\mathbf{W}=\{W_{t},t\geq 0\}$ is defined on a probability space $(\Omega ,~\mathbf{F}^{W}=\{F_{t}^{W}\}_{t\geq 0},~\mathbb{P})$, with  $\mathbf{F}^{W}$ being a filtration generated by $\mathbf{W}$. $b_{t}$ and $\sigma_{t}$ are measurable processes with respect to $\mathbf{F}^{Y}=\{F_{t}^{Y},~0\leq t\leq T\}$, where $F_{t}^{Y}={\mathcal {F}} \{Y_{s},0\leq s\leq t\}$ is the family of sigma-algebras generated by $\mathbf{Y}$.

The stochastic process driven by SDE in \eqref{Mainmodel} appears frequently in mathematical finance. For example, if $b_{t} = Y_{t}$ and $\sigma_{t}=\sigma$, then $Y$ becomes an Ornstein-Uhlenbeck process, which was introduced in finance in \cite{vasicek1977equilibrium}. The work was extended to more general type of SDE in \cite{CKLS}, known as the CKLS (Chan-Karolyi-Longstaff-Sanders) process, with $b_{t}=Y_{t}$ and $\sigma_{t}=\sigma Y_{t}^{\gamma}$, which is defined as
\begin{equation}
dY_{t}=\lambda (\alpha -Y_{t})dt+\sigma Y_{t}^{\gamma }dW_{t},\quad 0\leq t\leq T.
\label{CKLS}
\end{equation}
Note that from \eqref{CKLS}, if $\gamma=0$ and $\sigma_{t}=\sqrt{|Y_{t}|}$, then this leads to a CIR (Cox-Ingersoll-Ross) process in \cite{cox2005theory}. These stochastic processes were widely used for modelling the term structure of interest rates, and for volatility estimation in a finite time span; see \cite{phillips2009two} and \cite{laurent2020volatility}.

Several studies were made based on the SDE in \eqref{Mainmodel}. Theoretical results on the existence of strong and weak solutions were studied in the monographs in \cite{LipSh2001-1} and \cite{LipSh2001-2}, and statistical inferences for stochastic models described by SDE driven by Brownian motion have been presented in the monographs by \cite{Arato1982}, \cite{Prakasa}, \cite{Kuhler}, \cite{LipSh2001-1}, \cite{Bishwal}, \cite{Kutoyants2013}, \cite{Tanaka2017} and \cite{IacusYoshida2018}. The asymptotic properties of the estimators for discretised SDE-driven process in \eqref{Mainmodel} under various discretisation schemes can be found in \cite{kelly2004estimation}, while the asymptotic properties of sequential estimators for discretised O-U process have been studied in \cite{LaiSiegmund}, \cite{ABreton1989}, and \cite{novikovparameter}. In continuous time, the asymptotic properties of the maximum likelihood estimators (MLEs) for the CKLS process in \eqref{CKLS} were studied in \cite{Sorenson}, \cite{kelly2004estimation}, \cite{Cetin2013}, \cite{NovShir2014}, \cite{Higa2014}, \cite{Mish2022}. In addition, the computation and asymptotic properties of the pseudo moment estimators were studied in \cite{mota2018pseudo}, while a least squares approach has been studied in \cite{hu2010parameter} for a fractional Ornstein-Uhlenbeck processes, later extended to non-ergodic O--U processes in \cite{belfadli2011nonergodic}, \cite{el2016least} and \cite{el2017statistical}. 

In this paper, we study the bias and the mean-squared error (MSE) of the maximum likelihood estimators $\theta_{T}$ of $\theta = (\lambda, \alpha)$, which are defined as
\begin{equation*}
bias(\theta _{T})=\mathbb{E}_{\lambda ,\alpha }[(\theta _{T}-\theta )],\quad mse(\theta_{T})=\mathbb{E}_{\lambda ,\alpha }[(\theta _{T}-\theta )^{2}],
\end{equation*}
with the assumption that the sample size $T$ is fixed.

Studies on bias of estimators in autoregressive processes began from the classical study of the estimation of the correlation coefficient. The analytical calculation of the bias in the correlation estimator was first explored in \cite{marriott1954bias} and \cite{kendall1954note} in a stationary AR(1) model.  In particular, \cite{kendall1954note} highlighted the limitations of asymptotic expansions when $\rho$ is near unity, due to the skewness in the distribution of the estimator and the slow decay of higher-order terms. These concerns suggest that alternative approaches may be necessary to improve the accuracy of bias corrections. For a Gaussian stationary AR(1) process, the overview of asymptotics of least squares estimators for three parameters, including the bias corrections for the slope parameter was provided in \cite{orcutt1969first}. These studies have triggered further developments for bias analysis of various types of estimators, including the maximum likelihood estimators (see \cite{quenouille1956notes}, \cite{cordeiro1994bias}). In the case of continuous time, the asymptotic expansions for the bias and MSE of the maximum likelihood estimators analytically for the SDE in \eqref{Mainmodel}, in particular for O--U and CIR processes, have been studied in \cite{tang2009parameter}. The computational challenges arising from the dependence of the mean reversion parameter estimator on the initial value and discretisation frequency, as well as its slow convergence for small parameter values, have been examined in \cite{bao2017distribution}. 

We will derive an integral transform using the change of measure and Ito's formula to compute $bias(\theta_{T})$ and $mse(\theta_{T})$ for fixed $T$ under the assumptions made in \eqref{Mainmodel}. This approach allows us to derive expansions for cases as $\lambda \rightarrow 0$ or at $\lambda =\frac{u}{T}$, where $u$ is a constant. When referring to the study in \cite{novikovparameter}, the first-order asymptotic approximations for $bias(\theta _{T})$ and $mse(\theta _{T})$ could be misleading for moderate values of $T$ or small values of $\lambda$. The bias of the MLE can be used in the bias correction to reduce the mean-squared error of the MLE. 

This paper is organised as follows. In Section \ref{mle}, we formulate some general results studied in \cite{novikovparameter} that are valid within the framework of the mean-reverting process in \eqref{Mainmodel}. These results provides the basis for studying the fixed sample properties of the MLE for the Ornstein-Uhlenbeck (O-U) process. In Section \ref{ou_results}, we study the properties of MLEs to provide the analytical formulae and numerical illustrations for the bias and the MSE of the MLEs, and its modification using the moment generating functions (MGFs). The derivations of these MGFs are included in the Appendix. In Section \ref{numerical_illustration}, we provide a numerical study to illustrate the behaviour of bias and MSE for a finite sample size $T$. Section \ref{conc_ch4} concludes this paper.

\clearpage

\section{Maximum Likelihood Estimators}\label{mle}

In this section, we derive the expressions of the MLEs for the two parameters $\lambda$ and $\alpha$, and discuss their properties.

\subsection{Notations and Assumptions}

We assume $\mathbf{P}=\left\{\mathbb{P}_{\lambda,\alpha}(\alpha), A\in F_{T}^{Y}\right\}$ be the probability measure generated by the process $Y$ in \eqref{Mainmodel}, with $\mathbb{E}_{\lambda,\alpha}$ representing the expectation with respect to $\mathbb{P}_{\lambda,\alpha}$. In addition, assume that for any $\lambda \in \Lambda$ and $\alpha \in \mathcal{A}$, where $\Lambda$ and $\mathcal{A}$ are open intervals in $\mathbb{R}$, the following conditions hold:
\begin{align}
    &\mathbb{P}_{\lambda, \alpha}\left\{\int_{0}^{T} |b_t|\; dt < \infty \right\} = \mathbb{P}_{\lambda, \alpha}\left\{\int_{0}^{T} \sigma_{t}^2\; dt < \infty \right\} = 1,\\
    &\mathbb{P}_{\lambda, \alpha} \left\{\int_{0}^{T} \sigma_{t}^{-2}b_{t}^{2} \; dt < \infty \right\} = 1.
\end{align}
While these conditions hold, the measures $\mathbb{P}_{\lambda, \alpha}, \;\lambda\in\Lambda, \;\alpha \in\mathcal{A}$ are equivalent, and this indicates that for any parameters $\lambda\in\Lambda, \;\lambda_{0}\in\Lambda, \;\alpha \in \mathcal{A}$ and $\alpha_{0}\in\mathcal{A}$,
\begin{equation}
    \mathbb{P}_{\lambda, \alpha}(A) = \mathbb{E}_{\lambda_0, \alpha_0} \left[\mathbb{I}_{A}Z_{T}(\lambda, \alpha, \lambda_0, \alpha_0)\right],
\end{equation}
where $\mathbb{I}_{A}$ is an indicator function of an event $A\in F_{T}^{Y}$, and $\mathbf{Z} = \{Z_{t}(\lambda, \alpha, \lambda_0, \alpha_0),$
$\;0\leq t \leq T\}$ is a likelihood ratio process, also known as Radon-Nikodym derivative. 

Under assumptions from \cite{LipSh2001-1}, we have
\begin{align}
    Z_{t}&(\lambda, \alpha, \lambda_0, \alpha_0) \nonumber\\
    &= \exp \left\{\int_{0}^{t} \frac{\lambda(\alpha-b_{s})-\lambda_{0}(\alpha_{0}-b_{s})}{\sigma_{s}^{2}} \; dY_{s} - \frac{1}{2}\int_{0}^{t} \frac{\lambda^{2}(\alpha-b_{s})^2 - \lambda_{0}^2(\alpha_{0} - b_{s})^{2}}{\sigma_{s}^{2}} \; ds \right\}. \label{lkhd_z}
\end{align}
To simplify \eqref{lkhd_z}, we introduce the following notations.
\begin{align}
    \mathbf{S}(\alpha) &= \left\{S_{t}(\alpha) = \int_{0}^{t} \sigma_{s}^{-2}(\alpha-b_{s})\; dY_{s}, \; 0\leq t \leq T\right\},\label{notationS}\\
    \mathbf{Q}(\alpha) &= \left\{Q_{t}(\alpha) = \int_{0}^{t} \sigma_{s}^{-2} (\alpha-b_{s})^2\; ds, \; 0\leq t \leq T\right\},\label{notationQ}\\
    \mathbf{M}(\alpha) &= \left\{M_{t}(\alpha) = \int_{0}^{t} \sigma_{s}^{-1}(\alpha-b_{s})\; dW_{s}, \; 0\leq t \leq T\right\}.
\end{align}
To align with \eqref{Mainmodel},
\begin{align*}
    S_{t}(\alpha) &= \int_{0}^{t} \sigma_{s}^{-2}(\alpha-b_{s})\; dY_{s} = \lambda \int_{0}^{t} \sigma_{s}^{-2}(\alpha-b_{s})^2\; ds + \int_{0}^{t} \sigma_{s}^{-1}(\alpha-b_{s})\; dW_{s}\\
    &= \lambda Q_{t}(\alpha) + M_{t}(\alpha),
\end{align*}
where $Q_{t}(\alpha)$ is the quadratic characteristic of the local martingale $M_{t}(\alpha)$ with respect to the measure $\mathbb{P}_{\lambda, \alpha}$ and the filtration $\mathbf{F}^{Y}$. Hence, the likelihood ratio process $\mathbf{Z}$ can be written in the following form:
\begin{equation}
    Z_{t}(\lambda, \alpha, \lambda_{0}, \alpha_{0}) = \exp \left\{\lambda S_{t}(\alpha) - \lambda_{0}S_{t}(\alpha_{0}) - \frac{1}{2}(\lambda^{2} Q_{t}(\alpha) - \lambda_{0}^{2} Q_{t}(\alpha_{0}))\right\}.
\end{equation}

\subsection{Maximum Likelihood Estimators}
\textbf{(A) When $\alpha$ is known}, we first need to solve the following equation:
\begin{equation}
    \frac{\partial}{\partial \lambda} Z_{t}(\lambda, \alpha, \lambda_0, \alpha_0) = Z_{t}(\lambda, \alpha, \lambda_0, \alpha_0)(S_{t}(\alpha)-\lambda Q_{t}(\alpha)) = 0.
\end{equation}
and hence, the maximum likelihood estimator $\hat{\lambda}_{T}(\alpha)$ of $\lambda$ is
\begin{equation}
\hat{\lambda}_{T}(\alpha )=\frac{S_{T}(\alpha )}{Q_{T}(\alpha )}=\frac{\alpha \int_{0}^{T}\sigma _{s}^{-2}\;dY_{s}-\int_{0}^{T}\sigma_{s}^{-2}b_{s}\;dY_{s}}{\int_{0}^{t}\sigma _{s}^{-2}(\alpha -b_{s})^{2}\;ds}.
\label{MLE2}
\end{equation}
Note, that using \eqref{Mainmodel}, we may obtain the following useful representation:
\begin{equation}
    \hat{\lambda}_{T}(\alpha) = \lambda + \frac{S_{T}(\alpha)-\lambda Q_{T}(\alpha)}{Q_{T}(\alpha)} = \lambda + \frac{M_{T}(\alpha)}{Q_{T}(\alpha)}.\label{mlelambda}
\end{equation}

\textbf{(B) When $\lambda$ is known}, with $\lambda\neq0$, we need to solve the following equation:
\begin{equation}
    \frac{\partial}{\partial \alpha} Z_{t}(\lambda, \alpha, \lambda_0, \alpha_0) = Z_{t}(\lambda, \alpha, \lambda_0, \alpha_0)\lambda \left(\int_{0}^{t}\sigma_{s}^{-2}\;dY_{s} - \lambda\left(\int_{0}^{t} \sigma_{s}^{-2}(\alpha-b_{s})\;ds\right)\right),
\end{equation}
and hence, the maximum likelihood estimator $\hat{\alpha}_{T}(\lambda)$ of $\alpha$ is
\begin{equation}
    \hat{\alpha}_{T}(\lambda) = \frac{\int_{0}^{T} \sigma_{s}^{-2}\; dY_{s} + \lambda \int_{0}^{T} \sigma_{s}^{-2} b_{s} \; ds}{\lambda\int_{0}^{T}\sigma_{s}^{-2}\;ds}.
\end{equation}
For simpler representation, we denote
\begin{align}
    \mathbf{B}(\lambda) &= \left\{B_{t}(\lambda)=\int_{0}^{t} \sigma_{s}^{-2}\;dY_{s} + \lambda \int_{0}^{t} \sigma_{s}^{-2} b_{s}\;ds, 0\leq t \leq T\right\},\label{notationB}\\
    \mathbf{R} &= \left\{R_{t} = \int_{0}^{t} \sigma_{s}^{-2}\; ds, 0\leq t\leq T\right\},\label{notationR}\\
    \mathbf{N} &= \left\{N_{t}=\int_{0}^{t}\sigma_{s}^{-1}\; dW_{s}, 0\leq t\leq T\right\},
\end{align}
where $\mathbf{R}$ is the quadratic characteristic of the local martingale $\mathbf{N}$ with respect to the measure $\mathbb{P}_{\lambda, \alpha}$, and the filtration $\mathbf{F}^{Y}$. Using these notations, we can rewrite $\hat{\alpha}_{T}(\lambda)$ in the following form:
\begin{equation}
    \hat{\alpha}_{T}(\lambda) = \frac{B_{T}(\lambda)}{\lambda R_{T}}=\alpha + \frac{N_{T}}{\lambda R_{T}}.\label{mlealpha}
\end{equation}

\textbf{(C) When both $\lambda$ and $\alpha$ are unknown}, the MLEs can be found by solving the system of two equations with respect to $\lambda$ and $\alpha$,
\begin{equation}
    \begin{cases}
        \frac{\partial}{\partial \lambda} Z_{t}(\lambda, \alpha, \lambda_0, \alpha_0) = 0,\\
        \frac{\partial}{\partial \alpha} Z_{t}(\lambda, \alpha, \lambda_0, \alpha_0) = 0.
    \end{cases}
\end{equation}

The following relations between the MLEs $\hat{\alpha}_{T}$ and $\hat{\lambda}_{T}$ are obtained:
\begin{align}
    \hat{\alpha}_{T} &= \frac{\int_{0}^{T} \sigma_{s}^{-2}\; dY_{s} + \hat{\lambda}_{T}\int_{0}^{T} \sigma_{s}^{-2}b_{s}\;ds}{\hat{\lambda}_{T}R_{T}},\\
    \hat{\lambda}_{T} &= \frac{\hat{\alpha}_{T}\int_{0}^{T}\sigma_{s}^{-2}\;dY_{s}-\int_{0}^{t}\sigma_{s}^{-2}b_{s}\;dY_{s}}{Q_{T}(\hat{\alpha}_{T})}.
\end{align}

Now, note that in the SDE in \eqref{Mainmodel}, we have
\begin{equation}
    \frac{1}{T}\int_{0}^{T} b_{s}\;ds = \alpha + \frac{y-Y_{T}}{T}-\frac{1}{T}\int_{0}^{T}\sigma_{s}\;dW_{s}.
\end{equation}
Hence, as $T\rightarrow\infty$, if
\begin{equation}
    \frac{Y_{T}}{T}\overset{a.s.}{\longrightarrow }0 \text{ and } \frac{1}{T}\int_{0}^{T}\sigma_{s}\;dW_{s}\overset{a.s.}{\longrightarrow }0, \label{cond1}
\end{equation}
where the conditions in \eqref{cond1} hold when $\mathbf{Y}$ is ergodic, e.g. in the case of CKLS model in \eqref{CKLS}, then we obtain a strongly consistent estimator $\bar{\alpha}_{T}$,
\begin{equation}
    \bar{\alpha}_{T}= \frac{1}{T}\int_{0}^{T} b_{s}\; ds\overset{a.s.}{\longrightarrow }\alpha.\label{mlealpha_str}
\end{equation}
Using $\bar{\alpha}_{T}$, we may obtain the modified MLE of $\lambda$:
\begin{equation}
    \bar{\lambda}_{T} = \hat{\lambda}_{T}(\bar{\alpha}_{T})=\frac{\bar{\alpha}_{T}\int_{0}^{T}\sigma_{s}^{-2}\;dY_{s}-\int_{0}^{T}\sigma_{s}^{-2}b_{s}\;dY_{s}}{Q_{T}(\bar{\alpha}_{T})}.
\end{equation}
This MLE $\bar{\lambda}_{T}$ is connected to the previous MLE $\hat{\lambda}_{T}(\alpha)$ in \eqref{mlelambda} via the following relation:
\begin{equation}
    \bar{\lambda}_{T} = \hat{\lambda}_{T}(\alpha)\frac{Q_{T}(\alpha)}{Q_{T}(\bar{\alpha}_{T})}+(\bar{\alpha}_{T}-\alpha)\frac{\int_{0}^{T}\sigma_{s}^{-2}\;dY_{s}}{Q_{T}(\bar{\alpha}_{T})}.\label{lambdabar_lambdahat}
\end{equation}

If $\bar{\alpha}_{T}$ is a strongly consistent estimator, then as $T\rightarrow\infty$, it can be expected that
\begin{equation*}
    \frac{Q_{T}(\alpha)}{Q_{T}(\bar{\alpha}_{T})}\overset{a.s.}{\longrightarrow }1, \; \text{and } \bar{\lambda}_{T}-\hat{\lambda}_{T}(\alpha)\overset{a.s.}{\longrightarrow }0.
\end{equation*}
This will be shown in Section \ref{ou_results}, that in the case of an ergodic O--U process, the estimators $\bar{\alpha}_{T}$ and $\bar{\lambda}_{T}$ are asymptotically efficient.

\subsection{Properties of MLEs}\label{prop_mle}

Using the model in \eqref{Mainmodel}, we define Theorems 1 and 2 on general properties of MLEs, $\hat{\lambda}_{T}(\alpha)$ and $\hat{\alpha}_{T}(\lambda)$ in \eqref{mlelambda} and \eqref{mlealpha}, assuming one of the parameters is known. These theorems are reformulations of Theorem 1 from the study in \cite{novikovparameter} in the case of model \eqref{Mainmodel}.

%Denote $f(x)\in C^{n}(D)$ for a class of $n$-times continuously differentiable functions, $x\in D\in \mathbb{R}$.

\textbf{Theorem 1.} Assume $\alpha \in \mathcal{A}$.
\textit{
\begin{itemize}
    \item[a)] For $\lambda\in\Lambda$, let
    \begin{equation*}
        \mathbb{P}_{\lambda, \alpha}\left\{\lim_{T\rightarrow\infty} Q_{T}(\alpha)=\infty\right\}=1.
    \end{equation*}
    Then, 
    \begin{equation*}
        \mathbb{P}_{\lambda, \alpha}\left\{\lim_{T\rightarrow\infty} \hat{\lambda}_{T}(\alpha)=\lambda\right\}=1.
    \end{equation*}
    \item[b)] Let $p>0$ and $\sup_{\lambda\in\Lambda}\mathbb{E}_{\lambda, \alpha}[Q_{T}^{-r}(\alpha)]<\infty$ for some $r>p$. Then,
    \begin{equation*}
        \sup_{\lambda\in\Lambda}\mathbb{E}_{\lambda, \alpha}|\hat{\lambda}_{T}(\alpha)-\lambda|^{p}<\infty.
    \end{equation*}
    \item[c)] Let $\sup_{\lambda\in\Lambda} \mathbb{E}_{\lambda, \alpha}[Q_{T}^{-p}(\alpha)]<\infty$ for some $p>1$. Then,
    \begin{equation*}
       \mathbb{E}_{\lambda, \alpha}[Q_{T}^{-1}(\alpha)]\in C^{1}(\Lambda),\quad bias(\hat{\lambda}_{T}(\alpha))=\frac{\partial}{\partial\lambda}\mathbb{E}_{\lambda, \alpha}[Q_{T}^{-1}(\alpha)].
    \end{equation*}
    \item[d)] Let $\sup_{\lambda\in\Lambda} \mathbb{E}_{\lambda, \alpha}[Q_{T}^{-p}(\alpha)]<\infty$ for some $p>2$. Then,
    \begin{equation*}
        \mathbb{E}_{\lambda, \alpha}[Q_{T}^{-2}(\alpha)]\in C^{2}(\Lambda),\quad mse(\hat{\lambda}_{T}(\alpha))=\mathbb{E}_{\lambda, \alpha}[Q_{T}^{-1}(\alpha)]+\frac{\partial^2}{\partial\lambda^2}\mathbb{E}_{\lambda, \alpha}[Q_{T}^{-2}(\alpha)].
    \end{equation*}
    \item[e)] Let $\mathbb{E}_{\lambda, \alpha}[Q_{T}(\alpha)]\in(0,\infty)$ and $\sup_{\lambda\in\Lambda} \mathbb{E}_{\lambda, \alpha}[Q_{T}^{-p}(\alpha)]<\infty$ for some $p>1$. Then, for any $F_{T}^{Y}$-measurable integrable estimator $\hat{\lambda}_{T}$,
    \begin{equation}
        mse(\lambda_{T})\geq \frac{1}{\mathbb{E}_{\lambda, \alpha}[Q_{T}(\alpha)]} \left(1+\frac{\partial}{\partial\lambda}bias(\lambda_{T})\right)^2,\label{cr_lambda}
    \end{equation}
    \noindent which is the Cramer-Rao lower bound for $\lambda$.
    \item[f)] If there exists a constant $\delta>0$, such that $\mathbb{P}_{\lambda}\left\{Y_{T}(\alpha)\geq\delta\right\}=1$, then, for $\lambda\in\Lambda$,
    \begin{equation*}
        \mathbb{E}_{\lambda, \alpha}\left[\sqrt{Q_{T}(\alpha)}\right]<\infty \Rightarrow \mathbb{E}_{\lambda, \alpha}\left[\exp\left\{\frac{1}{2}|\hat{\lambda}_{T}|\right\}\right]<\infty. 
    \end{equation*}
\end{itemize}}

\textbf{Theorem 2.} Assume $\lambda$ is fixed, and $\lambda\neq0$.
\textit{\begin{itemize}
    \item[a)] For $\alpha\in\mathcal{A}$, let
    \begin{equation*}
        \mathbb{P}_{\lambda, \alpha}\left\{\lim_{T\rightarrow\infty} R_{T}=\infty\right\}=1.
    \end{equation*}
    Then, 
    \begin{equation*}
        \mathbb{P}_{\lambda, \alpha}\left\{\lim_{T\rightarrow\infty}\hat{\alpha}_{T}(\lambda)=\alpha\right\}=1.
    \end{equation*}
    \item[b)] Let $p>0$ and $\sup_{\alpha\in\mathcal{A}}\mathbb{E}_{\lambda, \alpha}[R_{T}^{-r}]<\infty$ for some $r>p$. Then,
    \begin{equation*}
        \sup_{\alpha\in\mathcal{A}}\mathbb{E}_{\lambda, \alpha}|\hat{\alpha}_{T}(\lambda)-\alpha|^{p}<\infty.
    \end{equation*}
    \item[c)] Let $\sup_{\alpha\in\mathcal{A}} \mathbb{E}_{\lambda, \alpha}[R_{T}^{-p}]<\infty$ for some $p>1$. Then,
    \begin{equation*}
        \mathbb{E}_{\lambda, \alpha}[R_{T}^{-1}]\in C^{1}(\mathcal{A}),\quad bias(\hat{\alpha}_{T}(\lambda))=\frac{1}{\lambda}\frac{\partial}{\partial\alpha}\mathbb{E}_{\lambda, \alpha}[R_{T}^{-1}].
    \end{equation*}
    \item[d)] Let $\sup_{\alpha\in\mathcal{A}} \mathbb{E}_{\lambda, \alpha}[R_{T}^{-p}]<\infty$ for some $p>2$. Then,
    \begin{equation*}
        \mathbb{E}_{\lambda, \alpha}[R_{T}^{-2}]\in C^{2}(\mathcal{A}),\quad mse(\hat{\alpha}_{T}(\lambda))=\frac{1}{\lambda}\mathbb{E}_{\lambda, \alpha}[R_{T}^{-1}]+\frac{1}{\lambda^{2}}\frac{\partial^{2}}{\partial\alpha^2}\mathbb{E}_{\lambda, \alpha}[R_{T}^{-2}].
    \end{equation*}
    \item[e)] Let $\mathbb{E}_{\lambda, \alpha}[R_{T}]\in(0,\infty)$ and $\sup_{\alpha\in\mathcal{A}}\mathbb{E}_{\lambda, \alpha}[R_{T}^{-p}]<\infty$ for some $p>1$. Then, for any $F_{T}^{Y}$-measurable estimator $\alpha_{T}$,
    \begin{equation*}
        mse(\alpha_{T})\geq \frac{1}{\lambda^2 \mathbb{E}_{\lambda, \alpha}[R_{T}]}\left(1+\frac{\partial}{\partial\alpha}bias(\alpha_{T})\right)^2
    \end{equation*}
\end{itemize}}

\textbf{Remarks}
\begin{itemize}
    \item[1.] The Cramer-Rao lower bound in Theorem 1(e) was considered in \cite{LipSh2001-2} under more general assumptions.
    \item[2.] For a case when the explicit formula for the Laplace transform $\mathbb{E}_{\lambda, \alpha}[\exp\{-\mu Q_{T}(\alpha)\}],$ $\;\mu \geq 0$, is known, then the negative moments $\mathbb{E}_{\lambda, \alpha}[Q_{T}^{-p}(\alpha)]$ can be obtained using the following formula:
    \begin{equation*}
        \mathbb{E}_{\lambda, \alpha}[Q_{T}^{-p}(\alpha)] = \frac{1}{\Gamma(p)}\int_{0}^{\infty} \mu^{p-1}\mathbb{E}_{\lambda, \alpha}[\exp\{-\mu Q_{T}(\alpha)\}]\;d\mu,
    \end{equation*}
    where $\Gamma(p)$ is a gamma function, defined as $\Gamma(p)=\int_{0}^{\infty} t^{p-1} e^{-t}\;dt.$
\end{itemize}

\clearpage

\section{The O-U Process and Results}\label{ou_results}

For this section, we consider a special case of the \eqref{Mainmodel}, when the process is defined by the SDE with $b_{t}=Y_{t}$ and $\sigma_{t}=\sigma$,
\begin{equation}
    dY_{t}=\lambda(\alpha-Y_{t})dt + \sigma dW_{t},\; Y_0=y,\; 0\leq t\leq T,\label{ou_simple}
\end{equation}
which is an O-U process. We further study the properties of the MLEs described in Section \ref{mle}, along with a brief numerical result based on this special case.

\subsection{General Properties of O-U Process}

To simplify notations in our work, we consider the case $\sigma=1$. However, the results for general $\sigma>0$ can be obtained by simple rescaling.
As $T\rightarrow \infty$, when $\lambda>0$, the process $Y$ is strongly ergodic, and
\begin{equation*}
    Y_{T}\overset{d}{\longrightarrow } Y_{\infty} \sim \mathcal{N}\left(\alpha, \frac{1}{2\lambda}\right).
\end{equation*}
A very well-known explicit solution for the SDE in \eqref{ou_simple} is given by
\begin{equation*}
    Y_{t}=\alpha+(y-\alpha)e^{-\lambda t} + e^{-\lambda t} \int_{0}^{t}e^{\lambda s}\;dW_{s},
\end{equation*}
and this indicates that $Y$ is a Gaussian process with the following moments:
\begin{align*}
    \mathbb{E}_{\lambda, \alpha}[Y_{t}]&=\alpha + (y-\alpha)e^{-\lambda t},\\
    Var_{\lambda, \alpha}[Y_{t}] &= \frac{1}{2\lambda}(1-e^{-2\lambda t}),\\
    \mathbb{E}_{\lambda, \alpha}[(Y_{t}-\alpha)^2] &= (y-\alpha)^2 e^{-2\lambda t}+\frac{1}{2\lambda}(1-e^{-2\lambda t}),\\
    Cov(Y_{T}, Y_{t})&= e^{-\lambda(T-t)}Cov\left(\int_{0}^{T}e^{\lambda s}\;dW_{s},\int_{0}^{t}e^{\lambda s}\;dW_{s}\right)=\frac{1}{2\lambda}e^{-\lambda(T-t)}(1-e^{-2\lambda t}).
\end{align*}
From \eqref{ou_simple}, the \textit{integrated O-U process} is defined as
\begin{equation}
    \mathbf{I}=\left\{I_{t}=\int_{0}^{t}Y_{s}\;ds,\;t\geq0\right\},\; I_{t}=\alpha t + \frac{1}{\lambda}(W_{t}-Y_{t}+y),\label{iOU}
\end{equation}
and the moments of the process $\mathbf{I}$ and $Cov(I_{t}, Y_{t})$ are
\begin{align*}
    \mathbb{E}_{\lambda, \alpha}[I_{t}]&=\alpha t + \frac{1}{\lambda}(y-\alpha-(y-\alpha)e^{-\lambda t})=\alpha t+\frac{y-\alpha}{\lambda}(1-e^{-\lambda t}),\\
    Var_{\lambda, \alpha}[I_{t}]&= \frac{1}{\lambda^{3}}\left[\frac{1}{2}(1-e^{-2\lambda t})-2(1-e^{-\lambda t})+\lambda t\right],\\
    Cov(I_{t}, Y_{t})&= \frac{1}{\lambda^2}(1-e^{-\lambda t})-\frac{1}{2\lambda^{2}}(1-e^{-2\lambda t}).
\end{align*}

\subsection{Properties of the MLE $\hat{\alpha}_{T}(\lambda)$ and $\bar{\alpha}_{T}$}

For the case in the process in \eqref{ou_simple}, by definitions in \eqref{notationB}, \eqref{notationR}, we obtain
\begin{align*}
    R_{T}&=\int_{0}^{T} ds = T, \text{ (since $\sigma=1$)}\\
    B_{T}(\lambda)&= \int_{0}^{T} dY_{s} + \lambda \int_{0}^{T} Y_{s}\;ds = Y_{T}-y+\lambda \int_{0}^{T} Y_{s}\;ds,
\end{align*}
and hence, the MLE $\hat{\alpha}_{T}(\lambda)$ in \eqref{mlelambda} is represented as:
\begin{equation}
    \hat{\alpha}_{T}(\lambda)=\frac{B_{T}(\lambda)}{\lambda R_{T}} = \frac{Y_{T}-y+\lambda\int_{0}^{T}Y_{s}\;ds}{\lambda T}=\alpha + \frac{W_{T}}{\lambda T}.\label{mlealpha_ou}
\end{equation}
This indicates that the bias and the mean-squared errors are
\begin{align}
    bias(\hat{\alpha}_{T}(\lambda)) &= \mathbb{E}_{\lambda, \alpha}[\hat{\alpha}_{T}(\lambda)-\alpha] = 0,\\
    mse(\hat{\alpha}_{T}(\lambda))&=\frac{1}{\lambda^2 T^2}Var[W_{T}]=\frac{1}{\lambda^2 T}.
\end{align}
As $\frac{W_{T}}{T}\overset{a.s.}{\longrightarrow }0$, we can easily obtain the result from \eqref{mlealpha_ou}, that $\hat{\alpha}_{T}(\lambda)\overset{a.s.}{\longrightarrow }\alpha$. In addition, by Theorem 2(e) in Section \ref{prop_mle} , the Cramer-Rao lower bound for any estimator $\alpha_{T}$ was

\begin{equation*}
    mse(\alpha_{T})\geq \frac{1}{\lambda^2\mathbb{E}_{\lambda, \alpha}[R_{T}]}\left(1+\frac{\partial}{\partial\alpha} bias(\alpha_{T})\right)^2,
\end{equation*}

and since $R_{T}=T$ and $bias(\alpha_{T})=0$, it can be shown that the MLE $\hat{\alpha}_{T}(\lambda)$ is the most efficient estimator among all other unbiased estimators.

Now, with the estimator $\bar{\alpha}_{T}$ from \eqref{mlealpha_str}, for the case of the O-U process, the estimator is expressed as
\begin{equation*}
    \bar{\alpha}_{T}=\frac{1}{T}\int_{0}^{T}Y_{s}\;ds.
\end{equation*}
Note that this estimator follows a normal distribution, and does not need the knowledge of $\lambda$. However, considering this estimator in terms of the integrated O-U process in \eqref{iOU} and its moments,
\begin{equation}
    \bar{\alpha}_{T}=\frac{1}{T}\int_{0}^{T}Y_{s}\;ds=\frac{I_{T}}{T}=\alpha+\frac{1}{\lambda T}(W_{T}-Y_{T}+y).
\end{equation}
and hence, the mean and the variance are
\begin{align}
    \mathbb{E}_{\lambda,\alpha}[\bar{\alpha}_{T}]&=\alpha+\frac{y-\alpha}{\lambda T}(1-e^{-\lambda T}),\\
    Var_{\lambda, \alpha}[\bar{\alpha}_{T}] &= \frac{1}{\lambda^{3}T^2}\left[\frac{1}{2}(1-e^{-2\lambda T})-2(1-e^{-\lambda T}) + \lambda T\right]\nonumber\\
    &= \frac{1}{\lambda^2T}\left[\frac{1}{2\lambda T}(1-e^{-2\lambda T})-\frac{2}{\lambda T}(1-e^{-\lambda T}) + 1\right]\\
    &= \frac{1}{\lambda^2 T}\left(1+O_{T}\left(\frac{1}{T}\right)\right),\; T\rightarrow \infty.
\end{align}
Thus, the estimator $\bar{\alpha}_{T}$ is asymptotically efficient amongst all other unbiased estimators. Also, since $Y_{T}\overset{d}{\longrightarrow }Y_{\infty}$ and $\frac{W_{T}}{T}\overset{a.s.}{\longrightarrow }0$ as $T\rightarrow\infty$, the estimator $\bar{\alpha}_{T}$ is a strongly consistent estimator of $\alpha$.

\subsection{Properties of the MLE $\hat{\lambda}_{T}(\alpha)$}\label{ssec:prop_lambda}

For the case of O-U process in \eqref{ou_simple} with $\sigma=1$, two integrals $Q_{T}(\alpha)$ and $S_{T}(\alpha)$ from \eqref{notationQ} and \eqref{notationS} are expressed as
\begin{align}
    Q_{T}(\alpha) &= \int_{0}^{T} (\alpha-Y_{s})^2\; ds = \int_{0}^{T} (Y_{s}-\alpha)^2\; ds,\label{notationQ_ou}\\
    S_{T}(\alpha) &= \int_{0}^{T}(\alpha-Y_{s})\;dY_{s} = \frac{1}{2}(T-(Y_{T}-\alpha)^2 + (y-\alpha)^2), \label{notationS_ou}
\end{align}
where $S_{T}(\alpha)$ is expressed after applying the Ito's formula. Using \eqref{notationQ_ou} and \eqref{notationS_ou} and \eqref{mlelambda}, the MLE $\hat{\lambda}_{T}(\alpha)$ is
\begin{equation}
    \hat{\lambda}_{T}(\alpha)=\frac{S_{T}(\alpha)}{Q_{T}(\alpha)}=\frac{T-Y_{T}^2+2\alpha Y_{T}-2\alpha y+y^2}{2Q_{T}(\alpha)}.
\end{equation}

Now, the distribution of $\hat{\lambda}_{T}(\alpha)$ in the case of the O-U process can be defined through the cumulative distribution function,
\begin{equation*}
    \mathbb{P}_{\lambda, \alpha}\left\{\hat{\lambda}_{T}(\alpha)<x\right\}=\mathbb{P}_{\lambda,\alpha}\left\{S_{T}(\alpha)-xQ_{T}(\alpha)<0\right\},
\end{equation*}
which leads to looking into the distribution of $\zeta(x)$, such that
\begin{equation*}
    \zeta(x)=S_{T}(\alpha)-xQ_{T}(\alpha)=T-Y_{T}^2+2\alpha Y_{T}-2\alpha y+y^2-2xQ_{T}(\alpha),\;y\in\mathbb{R}.
\end{equation*}

The moment generating function (MGF) of the random variable $\zeta(x)$ can be found through the joint moment-generating function of the random variables $\{Y_{T}, Q_{T}(\alpha)\}$. The explicit formula for the joint MGF is given by
\begin{equation}
    \psi = \left\{\psi_{\lambda,\alpha}(z_1, z_2, v, \mu)=\mathbb{E}_{\lambda, \alpha}\left[\exp\left\{z_{1}Y_{T}+z_{2}\int_{0}^{T}Y_{s}\;ds - vY_{T}^{2}-\mu Q_{T}(\alpha)\right\}\right]\right\}\label{mgf_psi}
\end{equation}
with $z_{1}\in\mathbb{R}, z_{2}\in\mathbb{R}, v\geq0, \mu\geq0, \lambda>0$.

By finding the MGF of the random variable $\zeta(x)$ through $\psi$, one can numerically find $\mathbb{P}_{\lambda, \alpha}\left\{\hat{\lambda}_{T}(\alpha)<x\right\}$ through a proper inversion formula, such as the Gaver-Stehfest algorithm in \cite{kuznetsov2013convergence}. Other methods include considering the analytical continuation of $\psi_{\lambda, \alpha}(z_{1}, z_{2}, v, \mu)$ to the region of complex values of the parameter $(z_{1}, z_{2}, v, \mu)$, and use techniques for numerical inversion of Fourier transform.

By Theorem 1 in Section \ref{prop_mle}, it is required to use the Laplace transform of $Q_{T}(\alpha)$ to find the bias and the mean-squared error of $\hat{\lambda}_{T}(\alpha)$. The Laplace transform is defined below.

\textbf{Proposition 1.} Let $Y$ be an O--U process with $\sigma=1$, $\lambda\in\mathbb{R}$, $\alpha\in\mathbb{R}$. Then for $\mu\geq0$,
\begin{align}
    &\psi_{\lambda,\alpha}(0,0,0,\mu) = \mathbb{E}_{\lambda, \alpha}\left[\exp\left\{-\mu \int_{0}^{T} (Y_{s}-\alpha)^2\;ds\right\}\right]\label{psi_000mu}\\
    &= \frac{1}{\left(1 + \frac{\lambda - \kappa}{2\kappa} (1 - e^{-2\kappa T}) \right)^{1/2}}\exp \left\{ \frac{(\lambda-\kappa) T}{2\kappa} + \frac{(y-\alpha)^2}{2} \left( \lambda - \kappa + \frac{2e^{-2\kappa T} }{1+ \frac{\lambda - \kappa}{2\kappa} (1 - e^{-2\kappa T})} \right) \right\}\nonumber,
\end{align}
where $\kappa=\sqrt{2\lambda+\mu^2}$. The above result was presented in \cite{N1973} with $\lambda\in\mathbb{R}$, including non-ergodic case as well, and the result has been used in \cite{Novikov1971Seq} and \cite{Nov1972MathNotes} for finding asymptotic expansions of bias and MSE of $\hat{\lambda}_{T}(\alpha)$ for two cases: $\lambda>0$ and $\lambda\rightarrow0$, see details in \cite{novikovparameter}.

Note that in the case when $\alpha=0$ and $\lambda\rightarrow0$, we have $Y=\mathbf{B}=\left\{B_{t}=W_{t}+y,\;t\geq0\right\}$, and through the former formulas presented, a version of the Cameron-Martin formula in \cite{CamMart1944} can be obtained, which is expressed as:
\begin{equation}
    \mathbb{E}\left[\exp{\left\{-\mu \int_{0}^{T} (W_{t}+y)^2\;ds\right\}}\right] = \frac{\exp\left\{-\sqrt{\mu/8}y^2 \tanh(\sqrt{2\mu}T)\right\}}{\sqrt{\cosh(\sqrt{2\mu}T})}, \; \mu\geq0.
\end{equation}
Now, using $\psi_{\lambda,\alpha}(0,0,0,\mu)$, all moments of $\hat{\lambda}_{T}(\alpha)$ can be calculated. Since
\begin{equation}
    \mathbb{E}_{\lambda,\alpha}\left[Q_{T}^{-p}(\alpha)\right]=\frac{1}{\Gamma(p)}\int_{0}^{\infty} \mu^{p-1}\psi_{\lambda,\alpha}(0,0,0,\mu)\;d\mu<\infty,\;p>0,
\end{equation}
and using Theorem 1 in Section \ref{prop_mle}, the bias and the MSE can be obtained.
\begin{align}
    bias(\hat{\lambda}_{T}(\alpha))&=\frac{\partial}{\partial\lambda}\mathbb{E}_{\lambda,\alpha}\left[Q_{T}^{-1}(\alpha)\right]=\int_{0}^{\infty}\frac{\partial}{\partial\lambda}\psi_{\lambda,\alpha}(0,0,0,\mu)\;d\mu,\label{biaslambda}\\
    mse(\hat{\lambda}_{T}(\alpha))&=\int_{0}^{\infty}\left(\psi_{\lambda,\alpha}(0,0,0,\mu)+\mu \frac{\partial^2}{\partial\lambda^2}\psi_{\lambda,\alpha}(0,0,0,\mu)\right)\;d\mu.\label{mselambda}
\end{align}
These calculations can be easily performed using mathematical computation programs such as \texttt{Wolfram Mathematica}. Alternatively, in general, the moments of $\hat{\lambda}_{T}(\alpha)$ can be calculated through the use of the distribution of $\hat{\lambda}_{T}(\alpha)$, which involves additional integrations, see \cite{Tanaka2017} for example.

$\hat{\lambda}_{T}(\alpha)$ is a strongly consistent estimator. The non-asymptotic result of the MSE is obtained through the Cramer-Rao lower bound for any estimator $\lambda_{T}$ in \eqref{cr_lambda}, where
\begin{align*}
    \mathbb{E}_{\lambda,\alpha}[Q_{T}(\alpha)] = \mathbb{E}_{\lambda,\alpha}\left[\int_{0}^{T}(Y_{s}-\alpha)^2\;ds\right]
    = \frac{1}{2\lambda}\left[T+\frac{(1-e^{-2\lambda T})}{2\lambda T}+(y-\alpha)^2(1-e^{-2\lambda T})\right],
\end{align*}
which can be obtained directly from \eqref{ou_simple}, or by finding $\frac{\partial}{\partial\mu} \psi_{\lambda,\alpha}(0,0,0,\mu)|_{\mu=0}$.

For the case when $T\rightarrow\infty$, 
\begin{equation*}
    \frac{Q_{T}(\alpha)}{T}=\frac{1}{T}\int_{0}^{T} (Y_{s}-\alpha)^2\;ds \overset{a.s.}{\longrightarrow} \mathbb{E}_{\lambda,\alpha}[(Y_\infty-\alpha)^2]=\frac{1}{2\lambda},
\end{equation*}
which leads to the following result:
\begin{equation*}
    \mathbb{E}_{\lambda,\alpha}[Q_{T}^{-1}(\alpha)]=\frac{1}{T}(\mathbb{E}_{\lambda,\alpha}[(Y_{\infty}-\alpha)^2])^{-1}(1+o_{T}(1))=\frac{2\lambda}{T}(1+o_{T}(1)).
\end{equation*}
Using Theorem 1(c) in Section \ref{prop_mle}, the bias is
\begin{equation}
    bias(\hat{\lambda}_{T}(\alpha))=\frac{\partial}{\partial\lambda}\mathbb{E}_{\lambda,\alpha}[Q_{T}^{-1}(\alpha)]=\frac{2}{T}(1+o_{T}(1)),\label{biaslambda_asymp}
\end{equation}
and also, since
\begin{equation*}
    \mathbb{E}_{\lambda,\alpha}[Q_{T}^{-2}(\alpha)]=\frac{1}{T^2}(\mathbb{E}_{\lambda,\alpha}(Y_\infty-\alpha)^2)^{-2}(1+o_{T}(1))=\frac{(2\lambda)^2}{T^2}(1+o_{T}(1)),
\end{equation*}
using Theorem 1(d) in Section \ref{prop_mle}, the MSE is
\begin{equation}
    mse(\hat{\lambda}_{T}(\alpha))=\mathbb{E}_{\lambda,\alpha}[Q_{T}^{-1}(\alpha)]+\frac{\partial^2}{\partial\lambda^2}\mathbb{E}_{\lambda,\alpha}[Q_{T}^{-2}(\alpha)]=\frac{2\lambda}{T}(1+o_{T}(1)).\label{mselambda_asymp}
\end{equation}
The asymptotic expansions for bias and MSE of $\hat{\lambda}_{T}(\alpha)$ in \eqref{biaslambda_asymp} and \eqref{mselambda_asymp} for cases when $\lambda\rightarrow0$ or $T\rightarrow\infty$ are of interest in the discipline of econometrics, see \cite{PP1991} for example. Using \eqref{psi_000mu}, it can be seen that for $T\rightarrow\infty$,
\begin{align}
   \lim_{\lambda\rightarrow0}\mathbb{E}_{\lambda,\alpha}[\hat{\lambda}_{T}(\alpha)-\lambda]&=\int_{0}^{\infty}\left[\frac{\partial}{\partial\lambda}\psi_{\lambda,\alpha}(0,0,0,\mu)|_{\lambda=0}\right]\;d\mu =\frac{C_0}{T^2}(1+o_{T}(1)),\\
    \lim_{\lambda\rightarrow0}\mathbb{E}_{\lambda,\alpha}[(\hat{\lambda}_{T}(\alpha)-\lambda)^2]&=\int_{0}^{\infty}\left[\psi_{\lambda,\alpha}(0,0,0,\mu)+\mu \frac{\partial^2}{\partial\lambda^2}\psi_{\lambda,\alpha}(0,0,0,\mu)|_{\lambda=0}\right]\;d\mu\nonumber\\
    &=\frac{C_1}{T^2}(1+o_{T}(1)),
\end{align}
where $C_0=1.7814...$ and $C_1=13.2857...$. Similar arguments can be applied to obtain the second order asymptotic expansions for $bias(\hat{\lambda}_{T}(\alpha))$ and $mse(\hat{\lambda}_{T}(\alpha))$ for the case when $\lambda=\frac{u}{T},\;T\rightarrow\infty$, where $u$ is a constant.

\subsection{Properties of the MLE $\bar{\lambda}_{T}=\hat{\lambda}_{T}(\bar{\alpha}_{T}).$}

Using \eqref{mlelambda}, the modified MLE estimator $\bar{\lambda}_{T}$ is expressed as
\begin{equation}
    \bar{\lambda}_{T} = \frac{T-Y_{T}^2-2\bar{\alpha}_{T}Y_{T}-2\int_{0}^{T}Y_{s}^2\;ds + y^{2}}{2Q_{T}(\bar{\alpha}_{T})},\; Q_{T}(\bar{\alpha}_{T}) = \int_{0}^{T}Y_{s}^2\;ds-(\bar{\alpha}_{T})^2T.
\end{equation}
To calculate the bias and the mean-squared error of $\bar{\lambda}_{T}$ numerically, the joint MGF of random vector $\{Y_{T}, \bar{\alpha}_{T}, (\bar{\alpha}_{T})^2, \int_{0}^{T} Y_{s}^{2}\;ds\}$ is used, i.e.)
\begin{equation}
    \boldsymbol{\Psi}=\left\{\Psi_{\lambda,\alpha}(z_{1}, z_{2}, \mu)=\mathbb{E}_{\lambda,\alpha}[\exp{z_{1}Y_{T}+z_{2}\bar{\alpha}_{T}-\mu Q_{T}(\bar{\alpha}_{T})}]\right\}, z_{1}\in\mathbb{R}, z_{2}\in\mathbb{R},\mu\geq0.
\end{equation}
For this MGF, we again adopt the tools from \texttt{Wolfram Mathematica}, as it is reasonably fast for calculating the bias and the MSE of $\bar{\lambda}_{T}$, see Appendix. The following formula can be used for calculating moments of $\bar{\lambda}_{T}$:
\begin{align}           
\mathbb{E}_{\lambda,\alpha}&\left[\frac{(\bar{\alpha}_{T})^{q}(Y_{T})^{j}}{(Q_{T}(\alpha))^{r}}\right]=\frac{1}{\Gamma(r)}\mathbb{E}_{\lambda,\alpha}\left[\int_{0}^{\infty}\mu^{r-1}(\bar{\alpha}_{T})^{q}(Y_{T})^{j}\exp\left\{-\mu Q_{T}(\bar{\alpha}_{T})\right\}\;d\mu\right]\nonumber\\
    &= \frac{1}{\Gamma(r)}\int_{0}^{\infty}\mu^{r-1}\frac{\partial^{q}}{\partial z_{2}^{q}}\frac{\partial^{j}}{\partial{z_{1}^{j}}}\boldsymbol{\Psi}\bigg|_{z_1=0,z_2=0}\;d\mu, r>0,\;q=0,1,\ldots,\;j=0,1,\ldots\label{mgf_Psi};
\end{align}
Through the formula in \eqref{mgf_Psi}, assuming $Y_{0}=0$, we have
\begin{equation*}
bias(\bar{\lambda}_{T})=\int_{0}^{\infty}\left(\frac{\partial}{\partial z_{2}}\frac{\partial}{\partial z_{2}}\boldsymbol{\Psi}-\frac{1}{2}\frac{\partial^2}{\partial z_{1}^2}\boldsymbol{\Psi}+\frac{T}{2}\boldsymbol{\Psi}\right)\bigg|_{z_1=0,z_2=0}\; d\mu -\lambda,
\end{equation*}
and similar formulas can be used to calculate $mse(\bar{\lambda}_{T})$ and higher moments of $\bar{\lambda}_{T}$.
Now, from \eqref{lambdabar_lambdahat}, since $\sigma=1$ and
\begin{equation*}
    \int_{0}^{T} \sigma_{s}^{-2}\;dY_{s}=\int_{0}^{T}\;dY_{s}=Y_{T}-y,
\end{equation*}
we get the following relation:
\begin{equation}
    \bar{\lambda}_{T}=\hat{\lambda}_{T}(\bar{\alpha}_{T})=\hat{\lambda}_{T}(\alpha)\frac{Q_{T}(\alpha)}{Q_{T}(\bar{\alpha}_{T})}+(\bar{\alpha}_{T}-\alpha)\frac{Y_{T}-y}{Q_{T}(\bar{\alpha}_{T})}, \label{modified_lambdabar-lambdahat}
\end{equation}
and through \eqref{modified_lambdabar-lambdahat}, it can be seen that $\bar{\lambda}_{T}$ is a strongly consistent estimator, since
\begin{equation*}
    \hat{\lambda}_{T}(\alpha)\overset{a.s.}{\longrightarrow}\lambda,\; \bar{\alpha}_{T}\overset{a.s.}{\longrightarrow }\alpha,\;Y_{T}\overset{d}{\longrightarrow }Y_{\infty}, \text{ as } T\rightarrow \infty,
\end{equation*}
and in the ergodic case, we get
\begin{equation*}
    \frac{Q_{T}(\alpha)}{Q_{T}(\bar{\alpha}_{T})}=\frac{\int_{0}^{T}(\alpha-Y_{s})^2\;ds}{\int_{0}^{T}(\bar{\alpha}_{T}-Y_{s})^2\;ds}\overset{a.s.}{\longrightarrow }1,\; \frac{Y_{T}-y}{\int_{0}^{T}(\bar{\alpha}_{T}-Y_{s})^2\;ds}\overset{a.s.}{\longrightarrow }0.
\end{equation*}
It can also be shown that as $T\rightarrow\infty$,
\begin{align*}
    bias(\bar{\lambda}_{T})&=bias(\hat{\lambda}_{T}(\alpha))(1+o_{T}(1))=\frac{2}{T}(1+o_{T}(1)),\\
    mse(\bar{\lambda}_{T})&=mse(\hat{\lambda}_{T}(\alpha))(1+o_{T}(1))=\frac{2\lambda}{T}(1+o_{T}(1)),
\end{align*}
and hence, $\bar{\lambda}_{T}$ is asymptotically the most efficient estimator among other unbiased estimators.

\clearpage

\section{Numerical Results}\label{numerical_illustration}

In this section, we provide some graphical illustrations for the bias and the mean-squared error of $\hat{\lambda}_T{}(\alpha)$ in \eqref{biaslambda} and \eqref{mselambda}. For ergodic cases, using the facts in \eqref{biaslambda_asymp} and \eqref{mselambda_asymp}, we present two following functions:
\begin{equation*}
    f_{1}(T)=T\times bias(\hat{\lambda}_{T}(\alpha)),\;f_{2}(T)=\frac{T}{\lambda}\times mse(\hat{\lambda}_{T}(\alpha)),
\end{equation*}
so that as $T\rightarrow\infty$, two functions approach to 2. We assume that $\sigma=1$ and $Y_0=1$.

\begin{figure}[h!]
    \centering
    \includegraphics[width=0.9\linewidth]{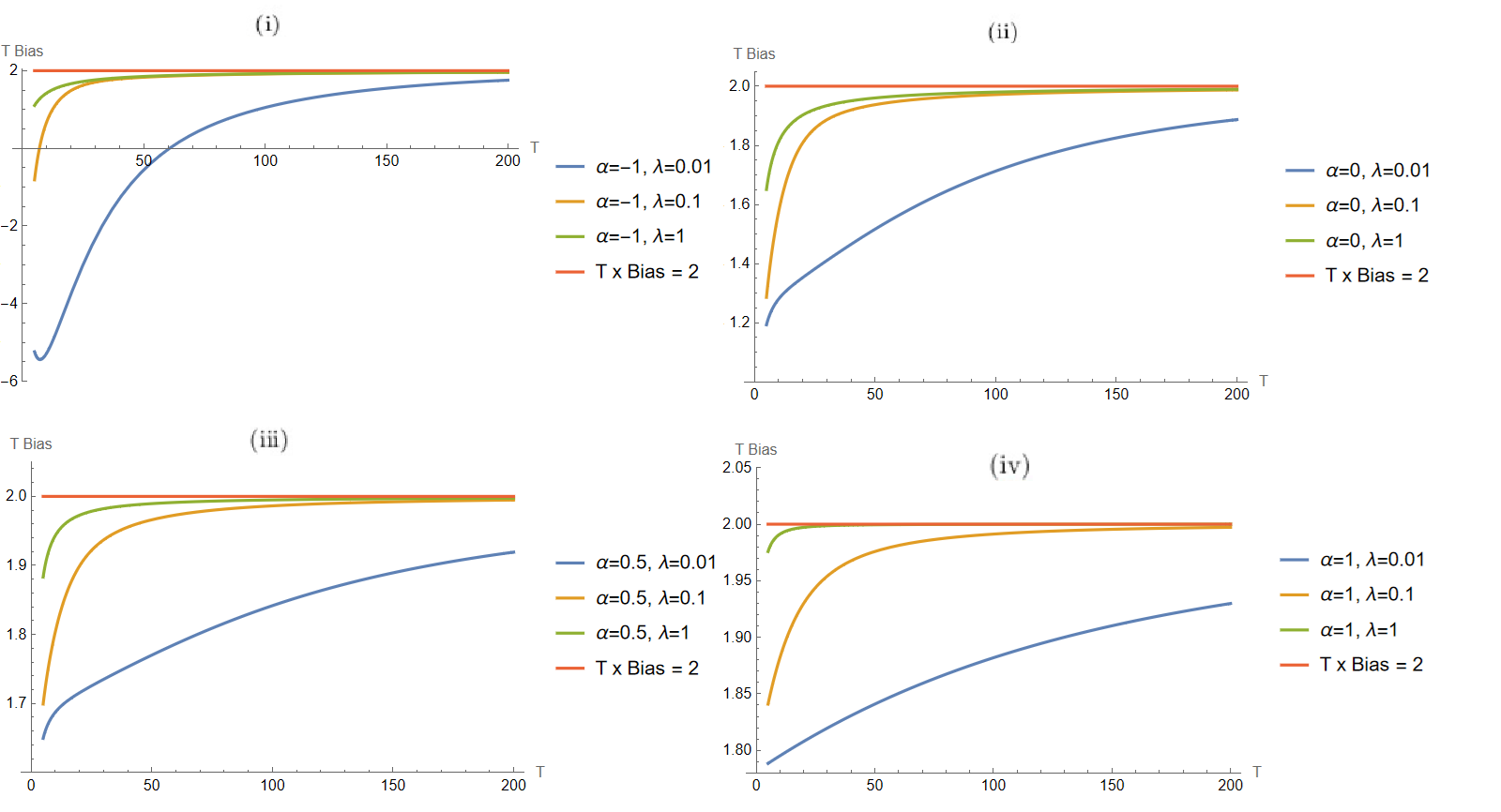}
    \caption{Plots of bias of $\hat{\lambda}_{T}(\alpha)$ for cases when (i) $\alpha=-1$, (ii)$\alpha=0$; (iii)$\alpha=0.5$; (iv) $\alpha=1$ with different values of $\lambda$}
    \label{fig:lambda_bias}
\end{figure}

\begin{figure}[h!]
    \centering
    \includegraphics[width=0.9\linewidth]{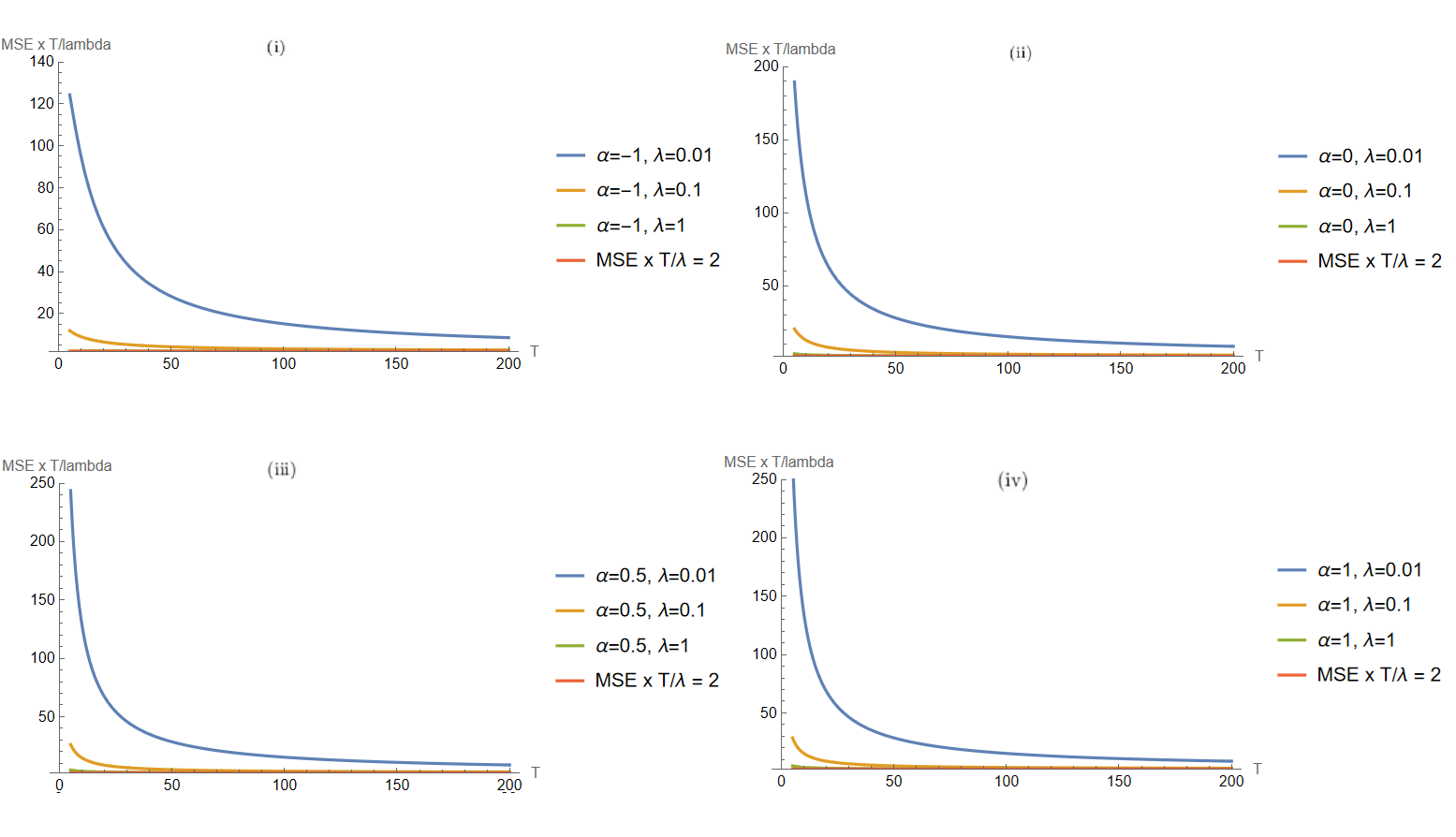}
    \caption{Plots of mean-squared errors of $\hat{\lambda}_{T}(\alpha)$ for cases when (i) $\alpha=-1$, (ii)$\alpha=0$; (iii)$\alpha=0.5$; (iv) $\alpha=1$ with different values of $\lambda$}
    \label{fig:mse_lambda}
\end{figure}

By \texttt{Wolfram Mathematica}, \eqref{biaslambda} and \eqref{mselambda} are used to approximate the bias and the mean-squared error numerically. Figures \ref{fig:lambda_bias} and \ref{fig:mse_lambda} depict bias and MSE in terms of two functions $f_{1}(T)$ and $f_{2}(T)$ across varying sample sizes $T$, and Table \ref{tab:bias_mse} in the Appendix provides a comparison of bias and MSE. Across all values of $\lambda$ and $\alpha$, the bias tends to decrease as $T$ increases, indicating better estimation with larger sample sizes. At smaller values of $\lambda$, particularly at $\alpha=1$, the MLE is underestimating $\lambda$, while for all other combinations, we see positive bias. The gap between bias values across different $\lambda$ and positive $\alpha$ shrinks, indicating that increasing $T$ reduces the dependence on $\lambda$. The bias decreases at a faster rate as both $T$ and $\lambda$ increase.

Looking into the MSEs, across all parameter settings, MSE decreases as $T$ increases, reflecting that larger sample sizes leads to more precise estimates. Larger $\lambda$ values tend to have higher MSE, particularly for small $T$, indicating that it leads to higher estimation error if the speed of mean-reversion is faster. At large $T$, the impact of both $\lambda$ and $\alpha$ on MSE is reduced, suggesting that sample size is the dominant factor in estimation accuracy. 


\begin{figure}[h!]
    \centering
    \includegraphics[width=\linewidth]{ 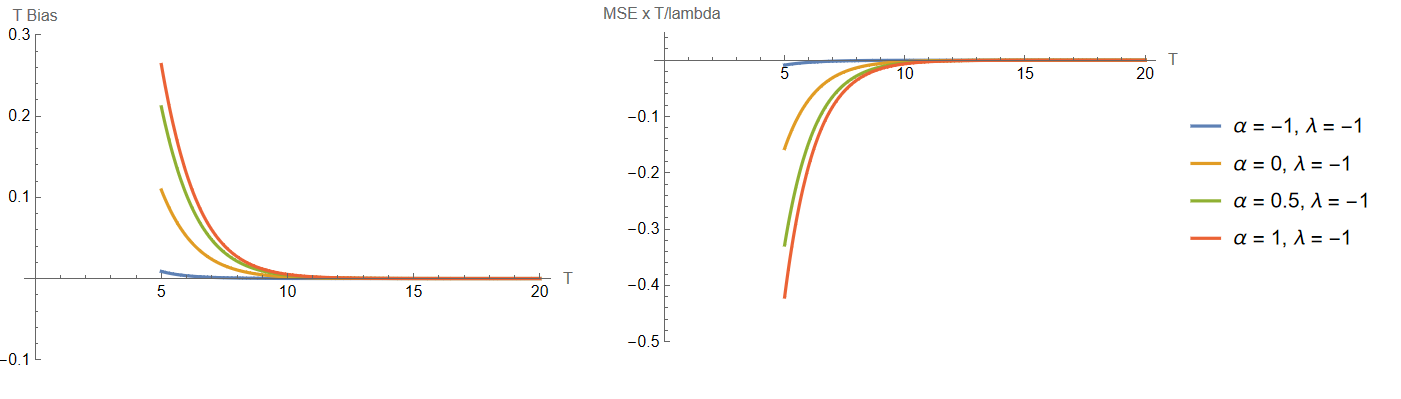}
    \caption{Plots of bias and mean-squared errors of $\hat{\lambda}_{T}(\alpha)$ for $\lambda=-1$ across different values of $\alpha$}
    \label{fig:bias_mse_nonerg}
\end{figure}

In contrast, the non-ergodic case displays qualitatively different behaviour. When $\lambda=-1$, we see that both $T\times bias(\hat{\lambda}_{T}(\alpha))$ and $\frac{T}{\lambda}\times MSE(\hat{\lambda}_{T}(\alpha))$ converge to zero, indicating a different scaling regime for the estimation error. Although the estimator no longer exhibits consistency in the classical ergodic sense, both scaled quantities decay rapidly as $T$ increases. As shown in Figure \ref{fig:bias_mse_nonerg}, the resulting curves are nearly indistinguishable for $T > 15$ across different values of $\alpha$, suggesting that the long-run mean has limited influence in this unstable regime.

\clearpage

\section{Conclusion}\label{conc_ch4}

This paper examined the theoretical properties of maximum likelihood estimators (MLEs) for a univariate, two-parameter mean-reverting process. Using the likelihood ratio process, we derived the expressions for MLEs of the parameters in cases when $\lambda$, $\alpha$, or both parameters were unknown, and formulated the general expressions for the bias and mean-squared errors (MSEs) of parameters.

From the general setting in \eqref{Mainmodel}, we simplified our model to follow an Ornstein-Uhlenbeck (O-U) process, with a constant volatility $\sigma_{t}=\sigma$. We derived the MLE of $\alpha$, the parameter representing the long-term mean, along with its bias and the MSE. For $\lambda$, the parameter representing the speed of mean-reversion, the bias and the MSE were derived based on the joint moment generating function in \eqref{mgf_psi}, and examined their asymptotic properties as the sample size $T$ increases.

In our numerical experiment, for an Ornstein-Uhlenbeck process with $\sigma=1, Y_0=1$, we provided graphical illustrations to demonstrate the dependency of bias and the MSE on the sample size and parameter values. Our findings indicate that as sample size increases, both bias and mean-squared error (MSE) consistently decrease across all parameter settings. While small values of $\lambda$ and $\alpha=-1$ led to a negative bias, and larger $\lambda$ led to higher MSE under small $T$, these effects diminish with larger samples. Figures \ref{fig:lambda_bias} and \ref{fig:mse_lambda} have shown that the bias and the MSE of $\hat{\lambda}_{T}(\alpha)$ approaches to $\frac{2}{T}$ and $\frac{2\lambda}{T}$, respectively, which verified the theoretical derivation in Section \ref{ssec:prop_lambda}. In contrast, under the non-ergodic setting $(\lambda=-1)$, the estimation error showed a different scaling behaviour, with bias and MSE converging to zero even at smaller sample sizes. Furthermore, as illustrated in Figure \ref{fig:bias_mse_nonerg}, the estimation results show minimal variation across different values of $\alpha$, indicating that the influence of the long-term mean becomes negligible in this unstable regime.

For further research, one could explore extending the current estimation framework and develop a discretised approximation method for mean-reverting diffusion processes under both ergodic and non-ergodic conditions. This includes incorporating fractional Ornstein-Uhlenbeck processes that allow for long-range dependence through fractional Brownian motion. Examining the effects of discretisation on the bias and MSE of parameter estimators in these general settings could enhance the understanding of inference under relaxed assumptions, which may contribute to the development of robust estimation procedures for models exhibiting structural complexities encountered in empirical applications.

\clearpage

%\setstretch{1.5}
%\small
\bibliographystyle{elsarticle-harv}
\bibliography{bibliography4}

\clearpage

\begin{appendix}
\section{Derivation of Moment Generating Function $\psi$}\label{append_4.A}

We derive the MGF of random vector $\{Y_{T}, Y_{T}^2, \int_{0}^{T} Y_{s}\;ds, \int_{0}^{T} Y_{s}^2\;ds\}$ in \eqref{mgf_psi} for the case when $b_{t}=Y_{t}$ and $\sigma_{t}=1$. The MGF was defined as
\begin{equation*}
    \psi=\left\{\psi_{\lambda,\alpha}(z_{1}, z_{2}, v, \mu)=\mathbb{E}_{\lambda,\alpha}\left[\exp\{z_{1}Y_{T}+z_{2}\int_{0}^{T}Y_{s}\;ds-vY_{T}^2-\mu Q_{T}(\alpha)\}\right]\right\},
\end{equation*}
where $z_{1}\in\mathbb{R}, z_{2}\in\mathbb{R}, v\geq0, \mu\geq0$.

\textbf{Theorem 3.} \textit{Let $\lambda > 0, Y_{0}=y$. Then,
\begin{align}
    \psi_{\lambda,\alpha}(z_{1}, z_{2}, v, &\mu)=D_{\lambda,\alpha}(T,y)\exp\left\{\frac{mr+m^2q+r^2d^2/2}{1-2qd^2}\right\}(1-2qd^2)^{-1/2},\label{psi_changem}\\
    \text{where }D_{\lambda,\alpha}&(T,y)= \exp\left\{(\lambda_0\alpha_0-\lambda\alpha)y+\frac{\lambda_0-\lambda}{2}y^2-\frac{T}{2}\left(\alpha^2\lambda^2+\lambda-\lambda_0^2\alpha_0^2-\lambda_0\right)\right\},\nonumber\\
    \lambda_0 &= \sqrt{\lambda^2+2\mu},\; \alpha_0=\frac{\lambda_0^2\alpha+z_{2}}{\lambda_0^2},\;r=z_{1}+\lambda\alpha-\lambda_0 \alpha_0,\;q=-\frac{\lambda_0-\lambda+2v}{2},\nonumber\\
    m&=\mathbb{E}_{\lambda_0,\alpha_0}[Y_{T}]=\alpha_0+(y-\alpha_0)e^{-\lambda t}, \;d^2=Var_{\lambda_0,\alpha_0}[Y_{T}]=\frac{1}{2\lambda_0}(1-e^{-\lambda_0 T}).\nonumber
\end{align}}

\begin{proof}
    The logarithm of the likelihood ratio process in \eqref{lkhd_z} is
    \begin{align*}
        &\log Z_{T}(\lambda,\alpha,\lambda_0,\alpha_0)\\
        &=\int_{0}^{T}\left[\lambda(\alpha-Y_{s})-\lambda_0(\alpha_0-Y_{s})\right]\;dY_{s}-\frac{1}{2}\int_{0}^{T}\left[\lambda^2(\alpha-Y_{s})^2-\lambda_{0}^{2}(\alpha_{0}-Y_{s})^2\right]\;ds\\
        &= \int_{0}^{T}(\lambda\alpha-\lambda_0\alpha_0)\;dY_{s}+\int_{0}^{T}(\lambda_0-\lambda)Y_{s}\;dY_{s}+\frac{\lambda_0^2 \alpha_0^2-\lambda^2\alpha^2}{2}T\\
        &+(\lambda^2\alpha-\lambda_0^2\alpha_0)\int_{0}^{T}Y_{s}\;ds
        + \frac{\lambda_0^2-\lambda^2}{2}\int_{0}^{T}Y_{s}^2\;ds\\ 
        &= \underbrace{(\lambda\alpha-\lambda_0\alpha_0)(Y_{T}-y)+\frac{\lambda_0-\lambda}{2}(T-Y_{T}^2+y^2)+\frac{\lambda_0^2\alpha_0^2-\lambda^2\alpha^2}{2}T}_\textrm{(i)}\\
        &+(\lambda^2\alpha-\lambda_0^2\alpha_0)\int_{0}^{T}Y_{s}\;ds+\frac{\lambda_0^2-\lambda^2}{2}\int_{0}^{T}Y_{s}^2\;ds.\;\left(\text{since} \; \int_{0}^{T} Y_{s}\;dY_{s}=\frac{T-Y_{T}^2+y^2}{2}\right)
    \end{align*}
By expanding (i), we get
\begin{align*}
    \text{(i)}&=\lambda\alpha Y_{T}-\lambda\alpha y-\lambda_0\alpha_0 Y_{T}+\lambda_0\alpha_0y+\frac{1}{2}\left(\lambda_0 T-\lambda T-\lambda_0 Y_{T}^2+\lambda Y_{T}^2+\lambda_0 y^2-\lambda y^2+\lambda_0^2\alpha_0^2 T-\lambda^2\alpha^2 T\right)\\
    &=\underbrace{(\lambda_0\alpha_0-\lambda\alpha)y+\frac{1}{2}(\lambda_0-\lambda)y^2-\frac{T}{2}\left(\lambda^2\alpha^2+\lambda-\lambda_0^2\alpha_0^2-\lambda_0\right)}_\textrm{$=\log D_{\lambda,\alpha}$}+(\lambda\alpha-\lambda_0\alpha_0)Y_{T}-\frac{1}{2}(\lambda_0-\lambda)Y_{T}^2
\end{align*}
Now, by the change of measure argument, since
$\frac{d\mathbb{P}_{\lambda,\alpha}}{d\mathbb{P}_{\lambda_0,\alpha_0}}=Z_{T}(\lambda,\alpha,\lambda_0,\alpha_0),$ we have $\mathbb{E}_{\lambda,\alpha}[g(\cdot)]=Z_{T}(\lambda,\alpha,\lambda_0,\alpha_0)\mathbb{E}_{\lambda_0,\alpha_0}[g(\cdot)].$

Hence, the MGF $\psi$ can be written as
\begin{align*}
    \psi_{\lambda,\alpha}&(z_1,z_2,v,\mu)=\mathbb{E}_{\lambda,\alpha}\left[\exp\left\{z_{1}Y_{T}+z_{2}\int_{0}^{T}Y_{s}\;ds-vY_{T}^2-\mu\int_{0}^{T}Y_{s}^2\;ds\right\}\right]\\
    &=Z_{T}(\lambda,\alpha,\lambda_0,\alpha_0)\mathbb{E}_{\lambda_0,\alpha_0}\left[\exp\left\{z_{1}Y_{T}+z_{2}\int_{0}^{T}Y_{s}\;ds-vY_{T}^2-\mu\int_{0}^{T}Y_{s}^2\;ds\right\} \right]\\
    &= D_{\lambda,\alpha}\times\mathbb{E}_{\lambda_0,\alpha_0}\bigg[\exp\bigg\{(\lambda\alpha-\lambda_0\alpha_0+z_{1})Y_{T}-\frac{\lambda_0-\lambda+2v}{2}Y_{T}^2\\
    &+(\lambda^2\alpha-\lambda_0^2\alpha_0+z_{2})\int_{0}^{T}Y_{s}\;ds + \frac{\lambda_0^2-\lambda^2-2\mu}{2}\int_{0}^{T}Y_{s}^{2}\;ds\bigg\}\bigg].
\end{align*}

Now, since $\lambda_0$ and $\alpha_0$ are arbitrary, we choose these two values, such that random variables $\int_{0}^{T}Y_{s}^2\;ds$ and $\int_{0}^{T}Y_{s}\;ds$ become zeros. That is, we have
\begin{equation*}
    \frac{\lambda_0^2-\lambda^2-2\mu}{2}=0,\;\lambda^2\alpha-\lambda_0^2\alpha_0+z_{2}=0.
\end{equation*}
This indicates that 
\begin{equation*}
\lambda_0=\sqrt{\lambda^2+2\mu},\; \alpha_0=\frac{\lambda^2{\alpha}+z_{2}}{\lambda_0^2}.
\end{equation*}

This simplifies the MGF $\psi$ to
\begin{equation}
    \psi_{\lambda,\alpha}(z_{1},z_{2},v,\mu)=D\times\mathbb{E}_{\lambda_0,\alpha_0}\left[\exp\left\{(\lambda\alpha-\lambda_0\alpha_0+z_{1})Y_{T}-\frac{\lambda_0-\lambda+2v}{2}Y_{T}^2\right\}\right],
\end{equation}
where $Y_{T}\sim\mathcal{N}(m,d^2)$.

We set $r = \lambda\alpha-\lambda_0\alpha_0+z_{1}$ and $q=-\frac{\lambda_0-\lambda+2v}{2}$, and find the joint MGF of $\left\{Y_{T}, Y_{T}^2\right\}$,
\begin{align}
    \mathbb{E}[\exp\left\{rY_{T}+qY_{T}^2\right\}]=\int_{-\infty}^{\infty} \frac{1}{\sqrt{2\pi d^2}}\exp\left\{rY_{T}+qY_{T}^2-\frac{1}{2d^2}(Y_{T}-m)^2\right\}\;dY_{T}    
\end{align}
By expanding and simplifying the exponents,
\begin{align*}
    rY_{T}+&qY_{T}^2-\frac{Y_{T}^2}{2d^2}+\frac{mY_{T}}{d^2}-\frac{m^2}{2d^2}=\left(q-\frac{1}{2d^2}\right)Y_{T}^2+\left(r+\frac{m}{d^2}\right)Y_{T}-\frac{m^2}{2d^2}\\
    &=AY_{T}^2+BY_{T}-\frac{m^2}{2d^2}, \; \left(A=q-\frac{1}{2d^2},\;B=r+\frac{m}{d^2}\right)\\
    &= A\left(Y_{T}^2+\frac{B}{A}Y_{T}\right)-\frac{m^2}{2d^2}=A\left[\left(Y_{T}+\frac{B}{2A}\right)^2-\frac{B^2}{4A^2}\right]-\frac{m^2}{2d^2}\\
    &= A\left(Y_{T}+\frac{B}{2A}\right)^2-\frac{B^2}{4A}-\frac{m^2}{2d^2}.
\end{align*}
Hence the integral is
\begin{align*}
    \frac{1}{\sqrt{2\pi d^2}}\int_{-\infty}^{\infty}& \exp\left\{A\left(Y_{T}+\frac{B}{2A}\right)^2\right\}\;dY_{T}\;\times\;\exp\left\{-\frac{B^2}{4A}-\frac{m^2}{2d^2}\right\}\\
    &= \frac{1}{\sqrt{2\pi d^2}}\sqrt{\frac{\pi}{-A}}\exp\left\{-\frac{B^2}{4A}-\frac{m^2}{2d^2}\right\}=\frac{1}{\sqrt{-2Ad^2}}\exp\left\{-\frac{B^2}{4A}-\frac{m^2}{2d^2}\right\}\\
    &= \left(-2d^2\left(q-\frac{1}{2d^2}\right)\right)^{-1/2}\exp\left\{-\frac{\left(r+\frac{m}{d^2}\right)^2}{4\left(q-\frac{1}{2d^2}\right)}-\frac{m^2}{2d^2}\right\}\\
    &= (1-2qd^2)^{-1/2}\exp\left\{\frac{mr+m^2q+r^2d^2/2}{1-2qd^2}\right\}.\\
    \therefore \psi_{\lambda,\alpha}(z_{1}, &z_{2}, v, \mu)=D_{\lambda,\alpha}(T,y)(1-2qd^2)^{-1/2}\exp\left\{\frac{mr+m^2q+r^2d^2/2}{1-2qd^2}\right\},
\end{align*}
which is the expression in \eqref{psi_changem}, and this completes the proof.
\end{proof}

\section{Derivation of Moment Generating Function $\Psi$}\label{append_4.B}

The MGF $\Psi$ in \eqref{mgf_Psi} for calculation of the bias and mean-squared error of the modified MLE $\hat{\lambda}_{T}$ was given by:
\begin{equation*}
    \boldsymbol{\Psi}=\left\{\Psi_{\lambda,\alpha}(z_{1}, z_{2}, \mu)=\mathbb{E}_{\lambda,\alpha}\left[\exp\left\{z_{1}Y_{T}+z_{2}\bar{\alpha}_{T}-\mu Q_{T}(\bar{\alpha}_{T})\right\}\right]\right\},
\end{equation*}
where $z_{1}\in\mathbb{R},\; z_{2}\in\mathbb{R}, \mu\geq0,\; Q_{T}(\bar{\alpha}_{T})=\int_{0}^{T}Y_{S}^2\;ds-(\bar{\alpha}_{T})^2 T.$

Firstly, we consider an auxiliary function,
\begin{equation*}
    \mathbf{U}=\left\{U_{\lambda,\alpha}(u_1,u_2,u_3,\mu)=\mathbb{E}[\exp\left\{u_1\xi_1+u_2\xi_2-\mu(\xi_1^2-u_3\xi_2^2)\right\}]\right\},
\end{equation*}
where $u_1\in\mathbb{R}, u_2\in\mathbb{R}, \mu\geq0, u_3\in\mathbb{G}$. $(\xi_1,\xi_2)$ is a two-dimensional Gaussian random vector, where
\begin{equation*}
    \left(\begin{array}{c}
        \xi_1\\
        \xi_2 
    \end{array}\right)\sim\mathcal{N}\left(
    \left(\begin{array}{c}
         m_1\\
         m_2
    \end{array}\right), \left(\begin{array}{cc}
       d_1  & d_{12} \\
       d_{12}  & d_2
    \end{array}\right)\right),
\end{equation*}
and $\mathbb{G}$ is the range, such that $\mathbb{G}=\left\{u_3:Q_{\lambda,\alpha}(u_1,u_2,v,\mu)<\infty\right\}$.

The MGF function $Q_{\lambda,\alpha}(u_1,u_2,\mu,u_3)$ is evaluated using \texttt{Wolfram Mathematica}, with the following expressions. First defining the bivariate-normal distribution,
\begin{equation*}
    (\xi_1,\xi_2)\sim \text{\texttt{N=MultinormalDistribution}}[m_1,m_2,\{\{d_1,d_{12}\},\{d_{12},d_2\}],
\end{equation*}
we find the MGF using the following command:
\begin{align*}
    &Q_{\lambda,\alpha}(u_1,u_2,\mu,u_3)\\
    &=\text{\texttt{Integrate}}[\text{\texttt{exp}}\{u_1 x_1+u_2 x_2-\mu(x_2^2-u_3x_1^2)\}\text{\texttt{PDF}}[\text{\texttt{N}},\{x_1,x_2\}],\{x_1,-\infty,\infty\},\{x_2,-\infty,\infty\}].
\end{align*}

Now, setting $\xi_1=Y_{T}, \xi_2=\bar{\alpha}_{T}, \lambda_0=\sqrt{\lambda^2+2\mu}, \alpha_0=\frac{\lambda^2\alpha+z_2}{\lambda_0^2}$, along with $Y_0=0$ and
\begin{align*}
    m_1&=\mathbb{E}_{\lambda_0,\alpha_0}[Y_{T}]=\alpha_0(1-e^{-\lambda_0T}), \\ d_{1}&=Var_{\lambda_0,\alpha_0}[Y_{T}]=\frac{1}{2\lambda_0}(1-e^{-\lambda_0T}),\\
    m_2&=\mathbb{E}_{\lambda_0,\alpha_0}[\bar{\alpha}_{T}]=\alpha_0-\frac{\alpha_0}{\lambda_0 T}(1-e^{-\lambda_0 T}), \\
    d_2 &= Var_{\lambda_0,\alpha_0}[\bar{\alpha}_{T}]=\frac{1}{\lambda_0^3 T^2}\left[\frac{1}{2}(1-e^{-2\lambda_0 T})-2(1-e^{-\lambda_0 T})+\lambda_0 T\right],\\
    d_{12}&=Cov(Y_{T},\bar{\alpha}_{T})=\frac{1}{T}\int_{0}^{T}Cov(Y_{T},Y_{t})\;dt=\frac{1}{2T\lambda^2}(1-e^{-\lambda T})^2,\\
    u_1&=z_1, u_2=z_2T, u_3=\mu T.
\end{align*}

\textbf{Theorem 4.} \textit{Let $\lambda>0, \; Y_0=0$. Then,
\begin{equation*}
    \Psi_{\lambda,\alpha}(z_1,z_2,\mu)=D_{\lambda,\alpha}(T,y)Q_{\lambda,\alpha}(z_1+\lambda\alpha-\lambda_0\alpha_0),z_2T,\mu,\mu T),
\end{equation*}
where $D_{\lambda,\alpha}$ is defined in \eqref{psi_changem}.}

\begin{proof}
    Theorem 3 can be modified in this case, and hence,
    \begin{align*}
        \Psi_{\lambda,\alpha}(z_1,z_2,\mu)&=D_{\lambda,\alpha}(T,y)\mathbb{E}_{\lambda_0,\alpha_0}\left[\exp\{(z_1+\lambda\alpha-\lambda_0\alpha_0)Y_{T}+z_2\bar{\alpha}_{T}+(\lambda_0-\lambda)Y_{T}^2+\mu T(\bar{\alpha}_{T})^2\}\right]\\
        &= D_{\lambda,\alpha}(T,y)Q_{\lambda,\alpha}((z_{1}+\lambda\alpha-\lambda_0\alpha_0),z_{2}\bar{\alpha}_{T},(\lambda_0-\lambda)Y_{T}^2,\mu T),
    \end{align*}
    where $\lambda_0=\sqrt{\lambda^2+2\mu}, \;\alpha_0=\frac{\lambda_0^2\alpha+z_2}{\lambda_0^2}$, as previously defined.
\end{proof}

\clearpage
\section{Summary of Numerical Results as a Table}
\begin{table}[h!]
    \centering
    \setlength{\tabcolsep}{8pt}
    \caption{Summary of bias (MSE) of $\hat{\lambda}_{T}(\alpha)$ for $\lambda=(0.01, 0.1, 1)$ and $\alpha=(-1, 0, 0.5, 1)$ at $T=(50, 75, 100, 125, 150, 175, 200)$}
    \resizebox{0.50\textwidth}{!}{\begin{tabular}{c|c|cccc}
        \toprule
        \(\lambda\) & \(T\) & \(\alpha=-1\) & \(\alpha=0\) & \(\alpha=0.5\) & \(\alpha=1\) \\
        \midrule
        \multirow{14}{*}{0.01} & \multirow{2}{*}{50} & -0.011 & 0.030 & 0.035 & 0.037 \\
                               &  & (0.006) & (0.006) & (0.006) & (0.006) \\
                               & \multirow{2}{*}{75} & 0.007 & 0.022 & 0.024 & 0.025 \\
                               &  & (0.003) & (0.003) & (0.003) & (0.003) \\
                               & \multirow{2}{*}{100} & 0.011 & 0.017 & 0.018 & 0.019 \\
                              &  & (0.002) & (0.002) & (0.002) & (0.002) \\
                               & \multirow{2}{*}{125} & 0.011 & 0.014 & 0.015 & 0.015 \\
                               &  & (0.001) & (0.001) & (0.001) & (0.001) \\
                               & \multirow{2}{*}{150} & 0.010 & 0.012 & 0.013 & 0.013 \\
                               &  & (0.001) & (0.001) & (0.001) & (0.001) \\
                               & \multirow{2}{*}{175} & 0.010 & 0.011 & 0.011 & 0.011 \\
                               &  & (0.001) & (0.001) & (0.001) & (0.001) \\
                               & \multirow{2}{*}{200} & 0.009 & 0.009 & 0.010 & 0.010 \\
                               &  & (0.0004) & (0.0004) & (0.0004) & (0.0004) \\
        \midrule
        \multirow{14}{*}{0.1} & \multirow{2}{*}{50} & 0.037 & 0.039 & 0.039 & 0.040 \\
                               &  & (0.008) & (0.009) & (0.009) & (0.009) \\
                               & \multirow{2}{*}{75} & 0.025 & 0.026 & 0.026 & 0.026 \\
                               &  & (0.005) & (0.005) & (0.005) & (0.005) \\
                               & \multirow{2}{*}{100} & 0.019 & 0.020 & 0.020 & 0.020 \\
                              &  & (0.003) & (0.003) & (0.003) & (0.003) \\
                               & \multirow{2}{*}{125} & 0.015 & 0.016 & 0.016 & 0.016 \\
                               &  & (0.002) & (0.002) & (0.002) & (0.002) \\
                               & \multirow{2}{*}{150} & 0.013 & 0.013 & 0.013 & 0.013 \\
                               &  & (0.002) & (0.002) & (0.002) & (0.002) \\
                               & \multirow{2}{*}{175} & 0.011 & 0.011 & 0.011 & 0.011 \\
                               &  & (0.02) & (0.02) & (0.02) & (0.02) \\
                               & \multirow{2}{*}{200} & 0.010 & 0.010 & 0.010 & 0.010 \\
                               &  & (0.001) & (0.001) & (0.001) & (0.001) \\
        \midrule
        \multirow{14}{*}{1}   & \multirow{2}{*}{50} & 0.037 & 0.039 & 0.040 & 0.040 \\
                               &  & (0.042) & (0.044) & (0.045) & (0.045) \\
                               & \multirow{2}{*}{75} & 0.025 & 0.026 & 0.027 & 0.027 \\
                               &  & (0.027) & (0.029) & (0.029) & (0.029) \\
                               & \multirow{2}{*}{100} & 0.019 & 0.020 & 0.020 & 0.020 \\
                              &  & (0.020) & (0.021) & (0.021) & (0.021) \\
                               & \multirow{2}{*}{125} & 0.016 & 0.016 & 0.016 & 0.016 \\
                               &  & (0.016) & (0.017) & (0.017) & (0.017) \\
                               & \multirow{2}{*}{150} & 0.013 & 0.013 & 0.013 & 0.013 \\
                               &  & (0.014) & (0.014) & (0.014) & (0.014) \\
                               & \multirow{2}{*}{175} & 0.011 & 0.011 & 0.011 & 0.011 \\
                               &  & (0.012) & (0.012) & (0.012) & (0.012) \\
                               & \multirow{2}{*}{200} & 0.010 & 0.010 & 0.010 & 0.010 \\
                               &  & (0.010) & (0.010) & (0.010) & (0.010) \\
        \bottomrule
    \end{tabular}}
    \label{tab:bias_mse}
\end{table}

\end{appendix}
    
\end{document}